\documentclass[11pt,reqno]{amsart}

\usepackage[papersize={8.5in,11in},width=6in,height=8in,centering]{geometry}
\usepackage{amsfonts}
\usepackage{amssymb}
\usepackage{amsthm}
\usepackage{amsmath}
\usepackage{graphicx}

\usepackage[all]{xy}

\theoremstyle{definition}
 \newtheorem{theorem}{\bf Theorem}[section]

 \newtheorem{corollary}[theorem]{Corollary}
 
\theoremstyle{definition}
 \newtheorem{example}[theorem]{Example}
 \newtheorem{remark}[theorem]{Remark}
 
 \newtheorem{proposition}[theorem]{Proposition}
\numberwithin{equation}{section}

\title[]{Concave symplectic toric fillings}

\author[A. Marinkovi\'c]{Aleksandra Marinkovi\'c} 

\email{aleksandra.marinkovic@matf.bg.ac.rs}

\begin{document}

\maketitle

\begin{abstract}
As shown by Etnyre and Honda (\cite{EH}), every contact 3-manifold admits infinitely many concave symplectic fillings that are mutually not symplectomorphic and not related by blow ups. In this note we refine this result in the toric setting by showing that every contact toric 3-manifold admits infinitely many concave symplectic toric fillings that are mutually not equivariantly symplectomorphic and not related by blow ups. The concave symplectic toric structure is constructed on certain linear and cyclic plumbings over spheres.
 \end{abstract}

\section{Introduction}
Let $(W,  \omega)$ be a symplectic manifold with a non-empty boundary $Y$ and let $V$ be a Liouville vector field  defined near the boundary and transversal to the boundary.  Then  $V$ induces a contact form $\lambda=i^{\ast}\iota_V\omega$ on $Y$, where $i:Y \to W$ is the inclusion map.
If $V$ points out of $Y$ then  the orientation on $Y$ given by $\lambda$ coincides with the boundary orientation on $Y$ induced from $W$ and 
$(W, \omega)$ is called a strong (convex) symplectic filling of $(Y, \ker \lambda)$.
If $V$ points inward, then 
$(W, \omega)$ is called a concave symplectic filling of $(Y, \ker \lambda)$. While the condition of being strongly symplectically fillable is quite restrictive (as it implies that the contact structure is tight \cite{Eli}),  every contact 3-manifold admits infinitely many concave symplectic fillings that are not related by blow ups (\cite[Theorem 1.3.]{EH}). 
These fillings  are constructed by Etnyre and Honda by gluing a Stein cobordism from a contact 3-manifold $Y$ to a certain Stein fillable 3-manifold $M$ and an infinity family of concave fillings of $M$. A Stein cobordism is built by attaching 2-handles along certain Legendrians to $Y\times[0,1]$, while concave fillings of $M$ are obtained by taking the complement of the family of Stein fillings  of $M$ symplectically embedded in  compact Kähler minimal surfaces. 

In this note we consider symplectic toric 4-manifolds with a non-empty boundary. Symplectic toric 4-manifolds are symplectic 4-manifolds equipped with an effective Hamiltonian $T^2=(\mathbb R/ \mathbb Z)^2 $ action. To every symplectic toric  manifold $(W, \omega)$ we associate a moment map  $H=(H_1,H_2),$ defined, up to a constant, by $\iota_{X_i}\omega=-dH_i,$ where $X_1$ and $X_2$ are the vector fields that  generate the toric action. 
Two symplectic toric manifolds are \emph{equivariantly symplectomorphic} if there exists a symplectomorphism between them that also preserves the toric actions.
 Suppose, a symplectic toric manifold $(W, \omega)$ admits a non-empty boundary $Y$ and a Liouville vector field transversal to the boundary. Then, by averaging, one obtains  a Liouville vector field that is invariant under the toric action on $W$. Thus, the induced contact form $\lambda$ on  $Y$ is also invariant under the toric action restricted to $Y$ and  $(Y, \xi)$ equipped with this action 
 is called a contact toric manifold.
 Moreover,  Cartan's formula $L_{X_i}\lambda=\iota_{X_i}d\lambda+d(\iota_{X_i}\lambda)$ leads to the unique moment map $H$ on $W$ that restricts to a moment map $H_{\lambda}=(\lambda(X_1), \lambda(X_2))$ on $Y.$ 
   
   The most common examples of symplectic toric manifolds with a contact toric boundary are toric domains, where a contact structure is necessarily convex. The first examples of symplectic toric manifolds with a concave contact toric boundary are provided by Nelson, Rechtman, Starkston, Tanny, Wang and the author in \cite[Theorem 4.1.]{MNRSTW}, where it is shown that any linear plumbing over spheres where at least one self-intersection number  is non-negative admits a symplectic toric structure with a concave contact toric boundary. 
   Further, in \cite[Theorem 1.4.]{MS} Starkston and the author  proved that every  contact 3-manifold with a non-free toric action admits  a concave symplectic toric filling by a certain linear plumbing over spheres.  In this note we extend this to an infinity family of such concave symplectic toric fillings and we also extend the result to contact 3-manifolds with a free toric action. The main result is the following.
   
   \begin{theorem} \label{thm1} Every contact toric 3-manifold admits infinitely many concave symplectic toric fillings that are mutually not equivariantly symplectomorphic and are not related by a sequence of  blow ups.
   \end{theorem}
      
\begin{figure}
\centering
\includegraphics[width=14cm]{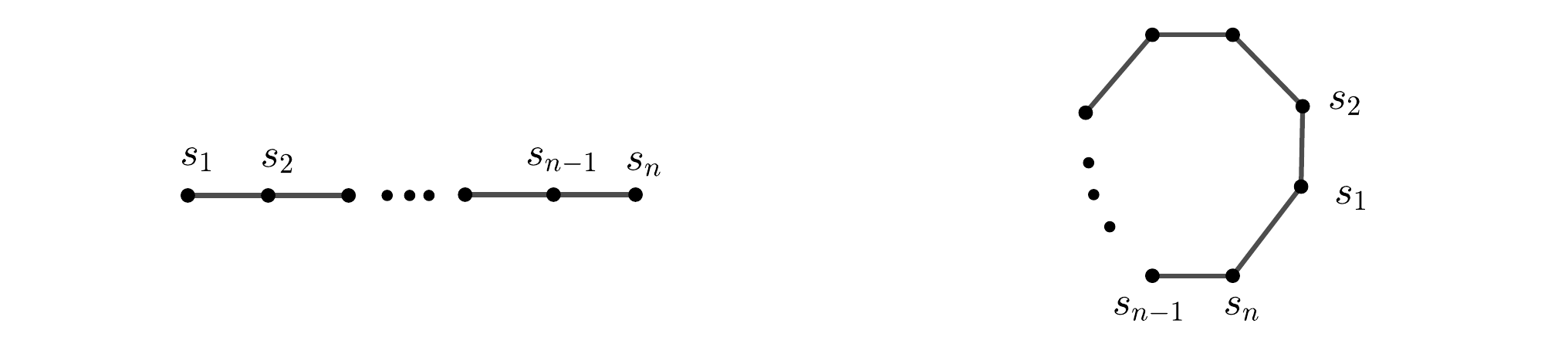}
\caption{Linear and cyclic plumbing graphs.}
\label{fig:graphs}
\end{figure}

The proof is divided in two parts. For a contact 3-manifold with a non-free toric action we use the combinatorial properties of the moment map image to find infinitely many distinct  linear plumbings over spheres that admit a concave symplectic toric structure with contactomorphic boundaries. As equivariant symplectomorphisms preserve the self-intersection numbers of base spheres it follows that distinct linear plumbings with constructed symplectic toric structures are not equivariantly symplectomorphic. For a contact 3-manifold with a free toric action we first provide conditions for a cyclic plumbing to admit a symplectic toric structure (Proposition \ref{thm0}) and then we find infinitely many desired cyclic plumbings over spheres, that are not equivariantly symplectomorphic and are not related by blow ups. See Figure \ref{fig:graphs} for the corresponding graphs,  where $s_1, \ldots, s_n  \in \mathbb Z$ denote the self-intersection numbers of the base spheres.


Finally, by ignoring the toric actions, we deduce the following corollary. 

   \begin{corollary}  There are three contact structures on any lens space $L(k, l)$, a universally tight one, half- and full-Lutz twist of the universally tight one,  that admit a concave symplectic filling by infinitely many distinct linear plumbings over spheres. All tight contact structures on $T^3$ admit a concave symplectic filling by infinitely many distinct cyclic plumbings over spheres.
   \end{corollary}
       
We remark that some of the constructed concave symplectic fillings may be symplectomorphic, but not by an equivariant symplectomorphism.
      
 \subsection*{Acknowledgments}
 The author is very grateful to  Laura Starkston for her encouragement, support, inspiration and for many useful comments. The author also thanks the anonymous referee  for  insightful comments and suggestions, which have improved the manuscript.
This research is partially supported by the Ministry of Education, Science and Technological Development, Republic of Serbia, through the project 451-03-136/2025-03/200104.

\section{Preliminaries on contact toric manifolds} \label{sec:prel}
  In this section we review some basic properties of contact toric manifolds. For more details we refer to   \cite{Lerman2}.
     
 A  contact manifold $(Y^{2n-1},\xi)$  equipped with an effective  $T^n=(\mathbb R/ \mathbb Z)^n $ action that preserves the contact structure $\xi$ is called \emph{a contact toric manifold}. 
To any invariant contact form $\alpha$ we associate a moment map  $H_{\alpha }=(H_1, \ldots, H_n): Y\to\mathbb R^{n}$  is uniquely
defined by 
$H_k=\alpha(X_k),$ $ k=1,\ldots, n,$
where $X_k$, $k=1,\ldots, n$ are the generators of the toric action. 
A moment map image is always  transversal to the radial rays emanating from the origin.  The union of the origin and the cone over a moment map image is called a moment cone and it depends only on the contact structure. By performing an automorphism of the torus $T^n$, the corresponding moment cone is changing by an $SL(n,\mathbb Z)$ transformation. If there exists a diffeomorphism between two contact toric manifolds that preserves the contact structures as well as the toric actions  we say that these contact toric manifolds are  \emph{equivariantly contactomorphic}.

     The toric action on a contact manifold may be free or non-free. We recall the classification in dimension 3, that will be relevant to prove the main results in the article.

      \begin{figure}
\centering
\includegraphics[width=9cm]{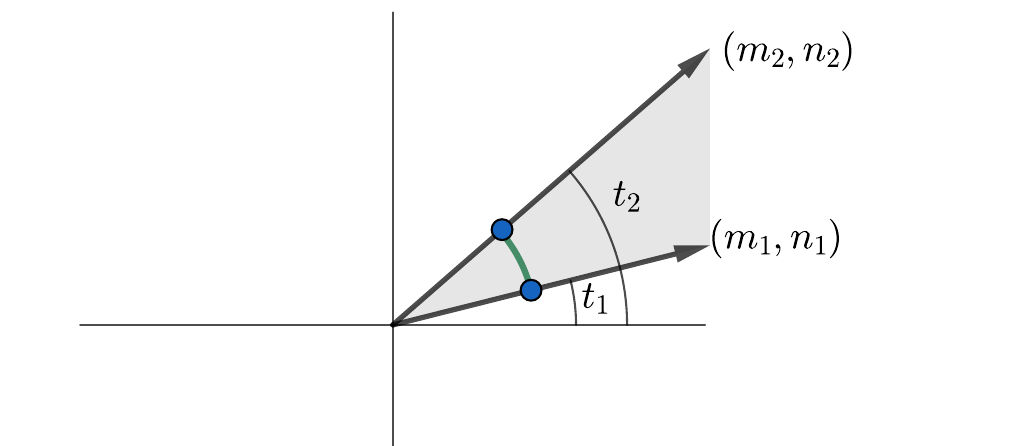}
\caption{A moment cone of a contact toric 3-manifold.}
\label{fig:cone}
\end{figure}

       \subsection{Contact 3-manifolds with a non-free toric action}    
     
\begin{theorem}(\cite[Theorem 2.18. (2)]{Lerman2})\label{theorem Lerman class}
Any  compact connected contact  manifold $(Y^3,\xi)$ with a non-free toric action is uniquely determined by two real numbers $t_1,t_2$ with $0\leq t_1<2\pi,t_1<t_2$ such that 
$\tan t_1$ and $\tan t_2$, when defined, are rational numbers.
\end{theorem}

 \begin{figure}[h!]
    \centering
    
        \begin{minipage}{0.75\textwidth}
        \centering
       \includegraphics[width=12cm]{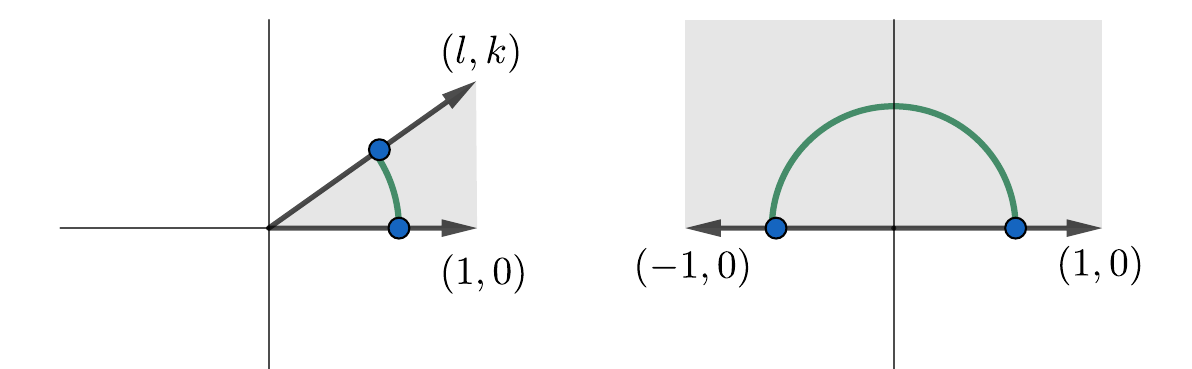}
\caption{Tight contact toric structure $  \xi_t$ on $L(k,l)$ (left) and $S^1\times S^2$ (right).}
        \label{fig:lens}
    \end{minipage}
      \hfill        
    \begin{minipage}{0.8\textwidth}
        \centering
        \includegraphics[width=\textwidth]{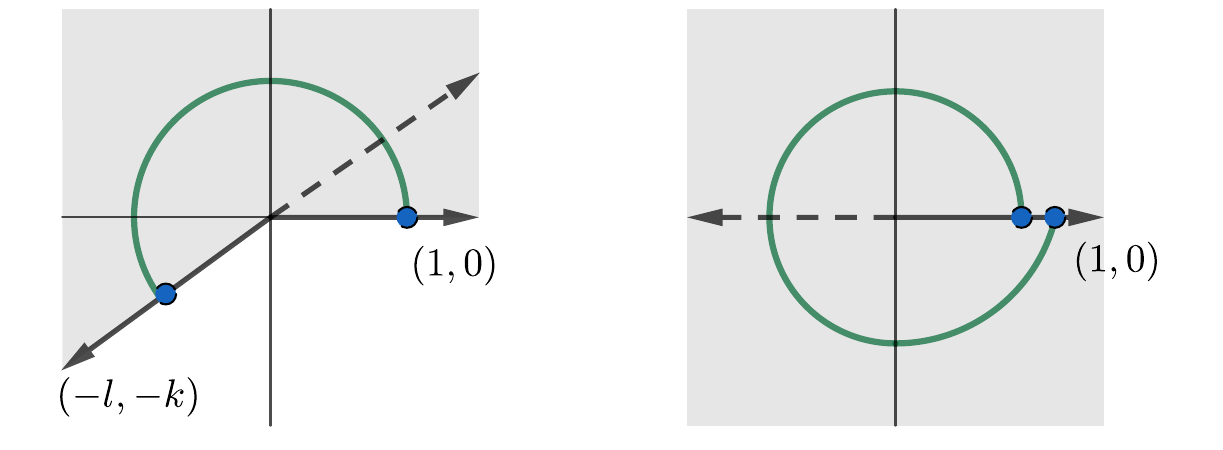}
        \caption{Overtwisted contact toric structure $  \xi_{ot1}$.}
        \label{fig:slika1}
    \end{minipage}
    \hfill
    \begin{minipage}{0.75\textwidth}
        \centering
        \includegraphics[width=\textwidth]{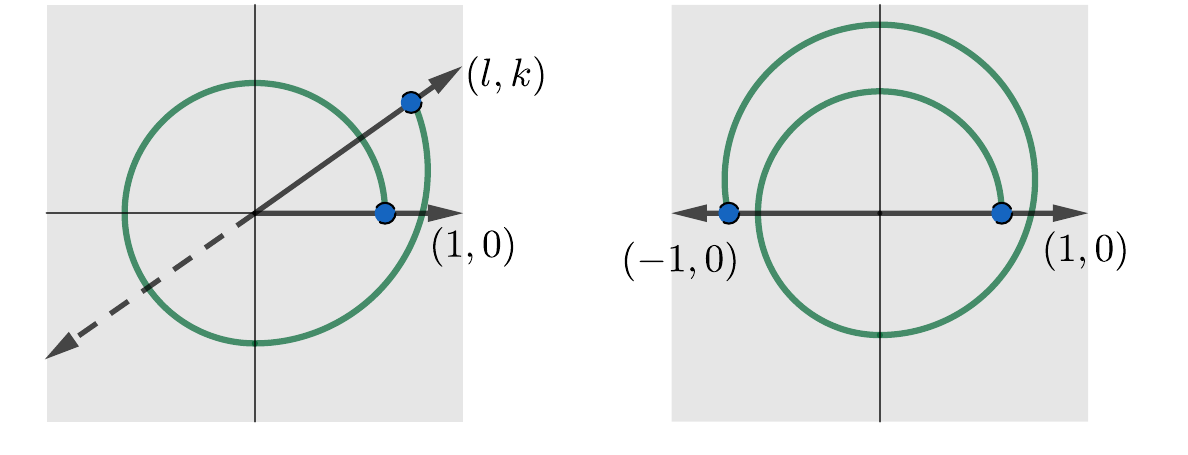}
        \caption{Overtwisted contact toric structure $  \xi_{ot2}$.}
        \label{fig:slika2}
    \end{minipage}

\end{figure} 

The numbers $t_1$ and $t_2$ are precisely the angles that the rays of the moment cone span with the positive part of the x-axis (see Figure  \ref{fig:cone}). According to Lerman, these rays are always given by rational slopes $(m_1,n_1)$ and $(m_2,n_2)$.
As $SL(2, \mathbb Z)$ transformations of moment cones preserve the corresponding contact toric structures, the numbers $t_1$ and $t_2$ are not unique numbers determining one contact toric structure.
In order to understand the contact toric structure start from $T^2\times [0,1]$ with coordinates $(\theta_1,\theta_2,t)$,  the contact structure $\ker(\cos (t_1(1-t)+t_2t)d\theta_1+\sin (t_1(1-t)+t_2t) d\theta_2)$ and the toric action given by the standard rotation of $T^2$ coordinates.
Collapse the  tori $T^2\times \{0\}$ and $T^2\times \{1\}$  along the circles of slopes $(-n_1,m_1)$ and $(n_2,-m_2)$, respectively. Note that these slopes are precisely the inward normal vectors to the rays that span the angles $t_1$ and $t_2$. The quotient space inherits a contact toric structure and $(Y, \xi)$ is equivariantly contactomorphic to it. Its moment  map image is depicted by green curve. Every point in the interior of the green curve corresponds to one $T^2$ orbit, while end points, depected by blue dots, correspond to circle orbits.

 Furthermore, the corresponding contact structure is tight  if and only if $t_2-t_1\leq \pi$  (see   \cite[Theorem 3.11.]{MNRSTW} ). 
 If $t_2-t_1< \pi$ the contact manifold is contactomorphic to a Lens space $L(k,l)$, for some $k>0, l \in \mathbb Z,$ with the unique universally tight contact structure $\xi_t$ (see Figure  \ref{fig:lens}), while
 if $t_2-t_1= \pi$ the contact manifold is contactomorphic to      $S^1\times S^2$ with the unique tight contact structures $\xi_t$. For more details see \cite[Theorem 1.1.]{MS}.
Next, if $\pi<t_2-t_1\leq 2\pi$ then,   according to \cite[Theorem 1.2]{MS}, the corresponding contact structure, denoted by $\xi_{ot1}$, is obtained by performing a half-Lutz twist to $\xi_t$.  
   By performing a sequence of full-Lutz twists to $\xi_{ot1}$ we obtain an isotopic overtwisted contact structure, with a different toric action. 
   The corresponding moment cone is defined by the same rays $(1,0)$ and
    $(-l,-k),$ however, the angle between the rays increases by $2n\pi ,$  where $n$ is the number of performed full-Lutz twists to $\xi_{ot1}$. 
   If $2\pi<t_2-t_1\leq 3\pi$ then the corresponding contact structure, denoted by $  \xi_{ot2}$  is obtained by performing a full-Lutz twist to $  \xi_t,$ or, by performing a half-Lutz twist to $ \xi_{ot1}.$
    By performing the sequence of full-Lutz twists to $  \xi_{ot2}$ we obtain isotopic overtwisted contact structures with different toric actions.
     The angle between the rays increases by $2n\pi ,$  where $n$ is the number of performed full-Lutz twists to $\xi_t$. 
As one full-Lutz twist coincides with two half-Lutz twists, we conclude the following, if $\pi<t_2-t_1$ the corresponding contact toric structure can be obtained from $\xi_t$ by performing the sequence of half-Lutz twists.

            \begin{figure}
\centering
\includegraphics[width=6.5cm]{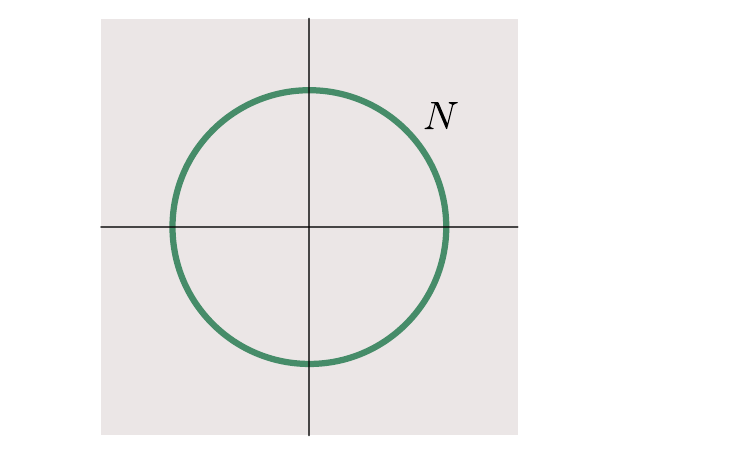}
\caption{A moment cone of $(T^3,\xi_N)$.}
\label{fig:free}
\end{figure}

   \subsection{Contact 3-manifolds with a free toric action}  
     
\begin{theorem}(\cite[Theorem 2.18. (1)]{Lerman2})\label{theorem Lerman class1}
Any  compact connected contact  manifold $(Y^3,\xi)$ with a free toric action is equivariantly contactomorphic to $(T^3,\xi_N=\ker(\cos N\theta d\theta_1+\sin N\theta d\theta_2))$, for some $N \in  \mathbb N, $ with the toric action given as the standard rotation of $(\theta_1,\theta_2)$ coordinates. 
\end{theorem}

     The corresponding moment map image  for all these contact structures is a unit circle and the moment cone is the whole space $\mathbb R^2.$ However, in contrast to the non-free toric action, where there are always two circle orbits, in the free case, the pre-image  of any point in the unit circle corresponds to $N$ tori $T^2$ in $T^3$ (see Figure \ref{fig:free}). Note also that all contact structures $\xi_N$, $N \in \mathbb N,$ are tight, while in the non-free case the moment cone 
      $\mathbb R^2$ corresponds only to overtwisted contact structures.

    \section{A symplectic toric structure on  linear and cyclic plumbings} \label{sec:plumb}
    Every plumbing is uniquely defined by
 a plumbing graph $\Gamma.$
      Each vertex of $\Gamma$ corresponds to a disc bundle over a surface of genus $g_i$ with self-intersection number $s_i$.  If two vertices are joined by an edge, then the corresponding two disc bundles are plumbed together. See Figure    \ref{fig:graphs} for a linear and cyclic graphs where all base surfaces are spheres. The cyclic plumbing graph is obtained from a linear graph by connecting the end vertices of a linear graph with an edge.	  To plumb two disc bundles, take a small discs in each base $D_1$ and $D_2$ and corresponding bundles $D_i(x)\times D^2(y)$. Glue $D_1\times D^2$ to $D_2\times D^2$ by a map 
 that switches the factors, orientation preserving or reversing.
  Associated to a (general) plumbing graph $\Gamma$ with $n$ vertices there is an  intersection form  $Q_\Gamma=[a_{ij}]_{i,j=1,  \ldots, n }$ defined by
  
 $  \bullet$ $a_{ii} =s_i$,  for every $1  \leq i   \leq n$;
 
$  \bullet$ $a_{ij} =0$, $i  \neq j $, if there is no edge between the vertices $i$ and $j$;

$  \bullet$ $a_{ij} =   \pm 1$, $i   \neq j,$ if there is an edge between the vertices $i$ and $j$, with an orientation $  \pm 1$ assigned.

  For more details on plumbings we refer to     \cite{Neu}.

 A plumbing where all gluing maps preserve the orientation admits a natural symplectic structure. Note that this means $a_{ij} \in  \{0,1 \},$ for all $i,j \in  \{1,  \ldots, n \}$, i.e. all plumbing edges are positive. A symplectic structure is constructed in the following way. On each base surface take a volume form and then
	choose
a symplectic structure on each disk bundle such that the zero section is a symplectic surface, and the disk
bundle is its standard small symplectic neighborhood.	
Assuming the symplectic areas of $D_1$ and $D_2$ are the same as the symplectic area of the $D^2$ fibers, the plumbing gluing is a symplectomorphism.	
Further,  such plumbings admit a contact type boundary.
If  $Q_{\Gamma}$ is negative-definite (for instance, if  $s_i \leq-2$, for all $i=1,  \ldots, n,$), then  the corresponding plumbing admits a convex contact type boundary (\cite{GS}),
while,
if there exists $z\in \mathbb R_{<0}^n$ such that $-Q_{\Gamma} z = a$, for some
$a=(a_1,\dots, a_n)\in \mathbb R_{>0}^n$,
then  the corresponding plumbing admits a concave contact type boundary (\cite{LM}) .

 \begin{figure}
\centering
\includegraphics[width=12cm]{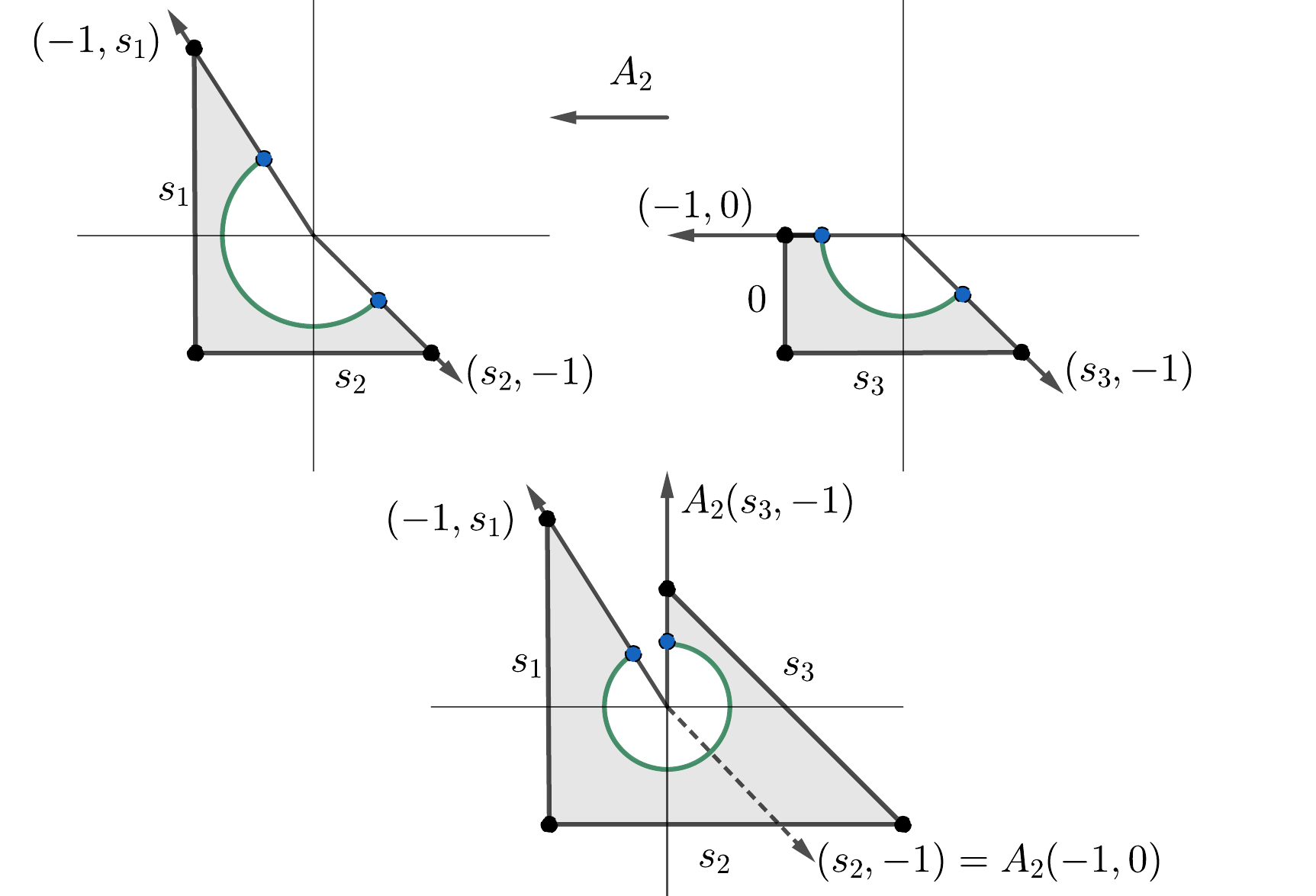}
\caption{A moment map image of a linear plumbing $(s_1,s_2,s_3)$ obtained by gluing the moment map images of linear plumbings $(s_1,s_2)$ and $(0,s_3)$.}
\label{fig:gluing}
\end{figure}   
Moreover, 
a linear plumbing $(s_1,  \ldots, s_n)$ over spheres where  $s_i \geq0$, for at least one index  $i=1, \ldots, n$,  admits a symplectic toric structure with a  concave contact toric boundary (\cite[Theorem 4.1.]{MNRSTW}).  
For a linear plumbing $(s_1,s_2)$ one constructs explicitly the moment map image, while for a linear plumbing $(s_1,  \ldots, s_n),$ $n>2,$ one performs the gluing of  moment map images of certain 2-vertices plumbings (see Figure  \ref{fig:gluing}). Namely, if $s_i\geq0$, then we consider the  plumbings 
$$ (s_1,0),\ldots, (s_{i-1},0),(s_i,s_{i+1}),(0,s_{i+2}),\ldots, (0,s_n),$$
and glue the corresponding moment map images starting from the moment map image of $(0,s_n) $ that will be glued to the left adjacent moment map image etc. Note that when we perform the gluing we keep
the segments along the rays that are glued to be the moment images of solid tori, as explained bellow in Remark \ref{remark:glue}.
 The gluing is uniquely given by the linear transformations of the moment map images $A_j=\begin{bmatrix}
 -s_j & -1 \\
1 &  0
\end{bmatrix}$, for all $j=2,  \ldots, n-1.$ Note that these linear maps  also induce the maps $A_j^{-T}$ on the toric fibers. 

 \begin{remark} \label{remark:equi} In general, every equivariant symplectomorphism of symplectic toric 4-manifolds is uniquely given by an $SL(2,\mathbb Z)$ transformation of the moment map image and every $SL(2,\mathbb Z)$ transformation defines an equivariant symplectomorphism.
  \end{remark}

In the  moment map image of a symplectic toric 4-manifold, the vertices  correspond to fixed points, the points on the interior of edges correspond to isotropic circle orbits and interior points correspond to Lagrangian 2-tori and the edges correspond to symplectic spheres. In the moment map image of a linear plumbing the edges correspond  precisely to the spheres in the base of the plumbing. The self-intersection number of the sphere corresponding to an edge $e$  is equal to the determinant of the inward normal vectors of the edges that are adjacent to $e$, taken in the clock-wise direction (\cite[Section 2.2.1.]{MNRSTW}). Since determinant is preserved under $SL(2, \mathbb Z)$ transformations, it follows that self-intersection numbers are preserved by equivariant symplectomorphisms.

 \begin{remark} \label{remark:glue} Note that in the moment map image of the plumbing $(s_1,s_2)$ the interior points on the rays correspond to circle orbits and black dots correspond to fixed points. Thus, the segment on each ray that starts at the black and ends at the blue dot corresponds to a  disc in a contact toric manifold. However, when we  perform the gluing, we require that the interior points on   the rays that will be glued correspond to 2-tori  and the black dots on the rays correspond to circle orbit. Thus, when we perform the gluing,  the segments from black to blue dots on the gluing rays correspond to solid tori.  Moreover,  the bottom edges in the moment map image of $(s_1,s_2)$ and in the moment map image of $(0,s_3)$ that correspond to the spheres with self-intersection numbers $s_2$ and $0$ are glued together into one edge that corresponds to one sphere in $(s_1,s_2,s_3)$ with self-intersection number $s_2$.
  \end{remark}

In the moment map image of this symplectic toric structure one can also read off the corresponding contact toric structure on the boundary, depicted by the green curve. 
The rays of the corresponding moment cone are given by

\begin{equation}\label{eq:rays}
R_1=(-1,s_1)\hskip2mm \textrm{and} \hskip2mm R_2=\begin{cases}
        (s_2,-1),&\mbox{ if $n=2$,  }\\  
  A_2  \cdots A_{n-1}(s_n,-1).&\mbox{ if $n\geq3.$}
 \end{cases}
\end{equation}

Note that, as explained in Section 2, the angle $t_2-t_1$ between the rays is very important, as it can change the contact structure.

       \subsection{Cyclic plumbings}    
We now extend the family of plumbings that admit a symplectic toric structure with a concave contact boundary. 

Note that  the linear plumbing $(s_1,  \ldots, s_n,s_{n+1}=0)$ admits a symplectic toric structure with a concave contact toric boundary, according to \cite[Theorem 4.1.]{MNRSTW}. Choose one such structure.
   
        \begin{proposition}\label{thm0}  Suppose that the rays of the contact toric boundary of the linear plumbing $(s_1,  \ldots, s_n,s_{n+1}=0)$ coincide and suppose that the  edges in the moment map image corresponding to $s_1$ and $s_{n+1}$ are parallel. 
     Then, the cyclic plumbing obtained from the  linear plumbing $(s_1,  \ldots, s_n)$ by plumbing the vertices $s_1$ and $s_n$ admits a structure of a symplectic toric manifold with a concave contact toric boundary such that the boundary is equivariantly contactomorphic to $(T^3,   \xi_N)$, for some $N\geq1$, with a free toric action. 
              \end{proposition}

         \begin{figure}
\centering
\includegraphics[width=13cm]{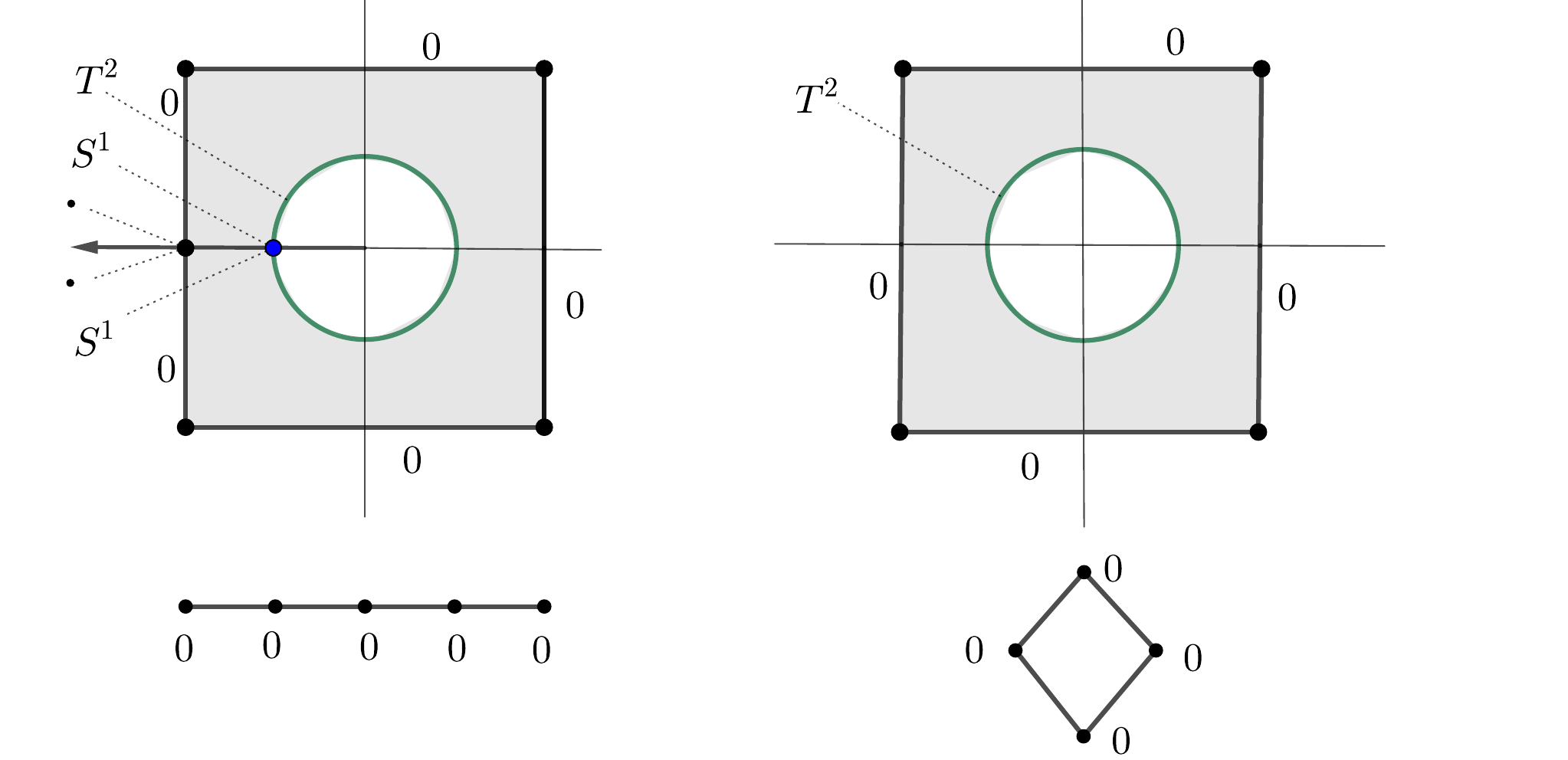}
\caption{A  symplectic toric structure on the cycling plumbing with $4$ vertices obtained from a linear plumbing with 5 vertices by gluing the moment map image along the rays. A concave contact boundary of the cyclic plumbing is equivariantly contactomorphic to $(T^3,   \xi_1)$.}
\label{fig:cyclic}
\end{figure}          

\begin{figure}
\centering
\includegraphics[width=12cm]{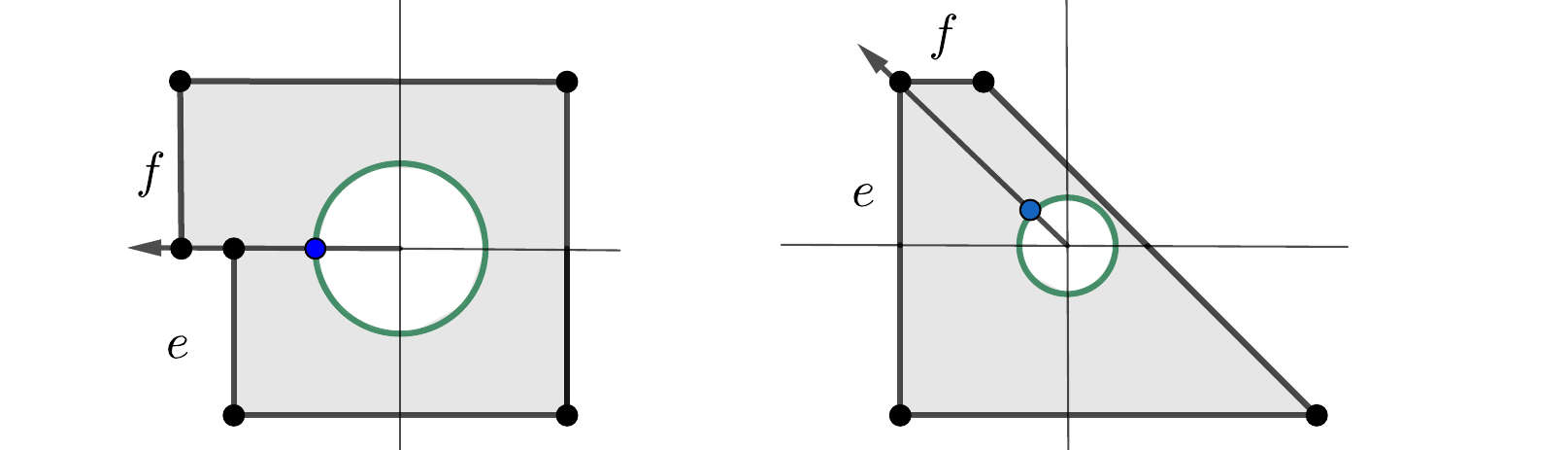}
\caption{The left linear plumbing  satisfies the conditions of Proposition \ref{thm0}, the right one doesn't. }
\label{fig:deform}
\end{figure}

           \begin{proof}                
            Since the rays  corresponding to the contact toric boundary of $(s_1,  \ldots, s_n,s_{n+1}=0)$ overlap,  we may suppose that the end vertices lying on the rays,  depicted by black dots, also overlap (see Figure \ref{fig:cyclic} on the left). Namely, this can be achieved by shrinking or extending the length of a certain edge in the moment map image (i.e. changing the symplectic area of the corresponding sphere).  
           See Figure \ref{fig:deform} on the left for a moment map image that can be deformed, by shrinking the edge adjacent to the edge $f$, to the moment map image shown on the left in Figure \ref{fig:cyclic}. 
           In terms of the  symplectic toric structures, shrinking and extending the edges, while preserving the inward normal vectors  and the rays of the moment cone, produce equivariantly deformation equivalent symplectic toric structures.

               As the end edges are parallel, the overlapping of the first and the last vertex implies that these edges lie on the same line and intersect in the overlapping vertex.
  Moreover, the condition that the rays of the contact toric boundary coincide implies that the numbers $t_1,t_2$ that determine the contact toric structure satisfy $t_2-t_1=2N  \pi,$ for some $N  \geq1.$

           Decompose the plumbing  $(s_1,  \ldots, s_n,0)$ into the sequence of plumbings
           $$ (s_1,0),\ldots, (s_{n-1},0), (s_n,0).$$ 
  We perform  the plumbings of the vertices $s_n$ and $s_1$ by gluing the moment map images of $(s_n,0)$ and $(s_1,0)$,  as in the case of two linear plumbings to obtain  a linear plumbing $(s_n,s_1,0)$. Namely, as explained in Remark  \ref{remark:glue}, two intermediate edges form one edge that corresponds to a sphere with a self-intersection number $s_1$. We then continue  as in the linear case by gluing the moment map image of $(s_1,0)$ to the moment map image of $(s_2,0)$ and obtain $(s_n,s_1,s_2,0)$. We continue inductively till we glue $(s_{n-1},0)$ to $(s_n,0)$. Therefore, obtained cyclic plumbing consists of $n$ vertices with self-intersection numbers $s_1,  \ldots, s_n.$

  Moreover, the blue dots on the boundary green curve in the moment map image of a linear plumbing will also correspond to 2-tori in the cyclic plumbing. Therefore, there are no circle orbits at the boundary of a cyclic plumbing. Since the moment map image of the boundary is a circle and all orbits are 2-tori the total space is diffeomorphic to $T^3.$ The number $N$ of full circles determines the contact structure $(T^3,   \xi_N )$.

       \end{proof}

         \begin{example}      
         The linear plumbing $(0,0,0,0,0)$ with a symplectic toric structure induces a symplectic toric structure on the cyclic plumbing $(0,0,0,0)$, see  Figure \ref{fig:cyclic}. 
  Since the angle between the rays is $2\pi$, the concave contact boundary of the linear plumbing is equivariantly contactomorphic to $(S^1 \times S^2,  \xi_{ot1})$,
  while the concave contact boundary of the cyclic plumbing is  equivariantly contactomorphic to $(T^3,   \xi_1 )$. In general, starting from a linear plumbing $(\underset{4N+1}{\underbrace{0, \ldots, 0}})$
  one obtains a cyclic plumbing with $4N$ vertices whose boundary is equivariantly contactomorphic to   $(T^3,   \xi_N )$.     
         \end{example}      
         
                  \begin{remark}      
The condition of parallel edges cannot be excluded from the statement of Proposition \ref{thm0}. Namely, the linear plumbing $(1,1,0,0)$
 admits a symplectic toric structure with a concave contact toric boundary, where the rays of the moment cone also overlap. However, the end edges, labeled by $e$ and $f$ in Figure \ref{fig:deform}, on the right are not parallel and the plumbing of the first and the last vertex cannot be performed in such a way to obtain a symplectic toric structure.

         \end{remark}              
\section{Proof of Theorem \ref{thm1} }  \label{sec:thm1}

\subsection{Contact  3-manifolds with a non-free toric action}
We first prove the theorem for  the tight contact toric structures and then for the overtwisted ones.

\noindent \textbf{1. Case} If $t_2-t_1=\pi$ the corresponding contact toric manifold is equivariantly contactomorphic to $S^1 \times S^2$ with the unique tight contact structure $\xi_t$ and the moment cone shown on the right in Figure \ref{fig:lens}.
By a simple computation that relies on (\ref{eq:rays}), we conclude that this  contact toric structure  can be realised as
a concave contact boundary of the  linear plumbings
         $$(n,0,-n), n  \geq0,$$
         with a symplectic toric structure shown on the left in Figure 9.   
           Since equivariant symplectomorphisms preserve self-intersection numbers of spheres corresponding to edges in the moment map image, it follows that these plumbings are mutually not equivariantly symplectomorphic.

 \begin{figure}[h!]     \label{fig:two}
     \centering
    
        \begin{minipage}{0.49\textwidth}
        \centering
       \includegraphics[width=\textwidth]{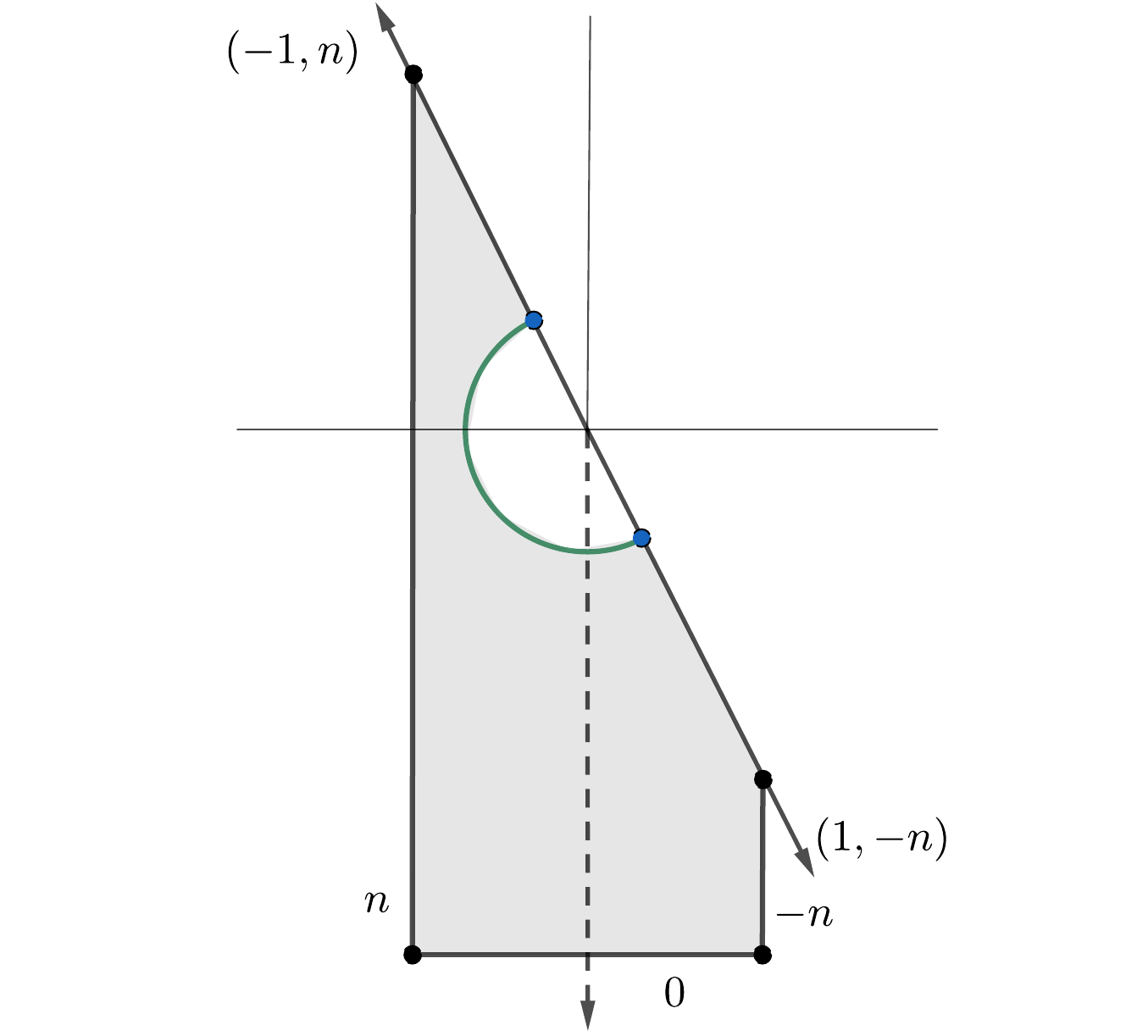}
    
    \end{minipage}
      \hfill        
    \begin{minipage}{0.49\textwidth}
        \centering
        \includegraphics[width=\textwidth]{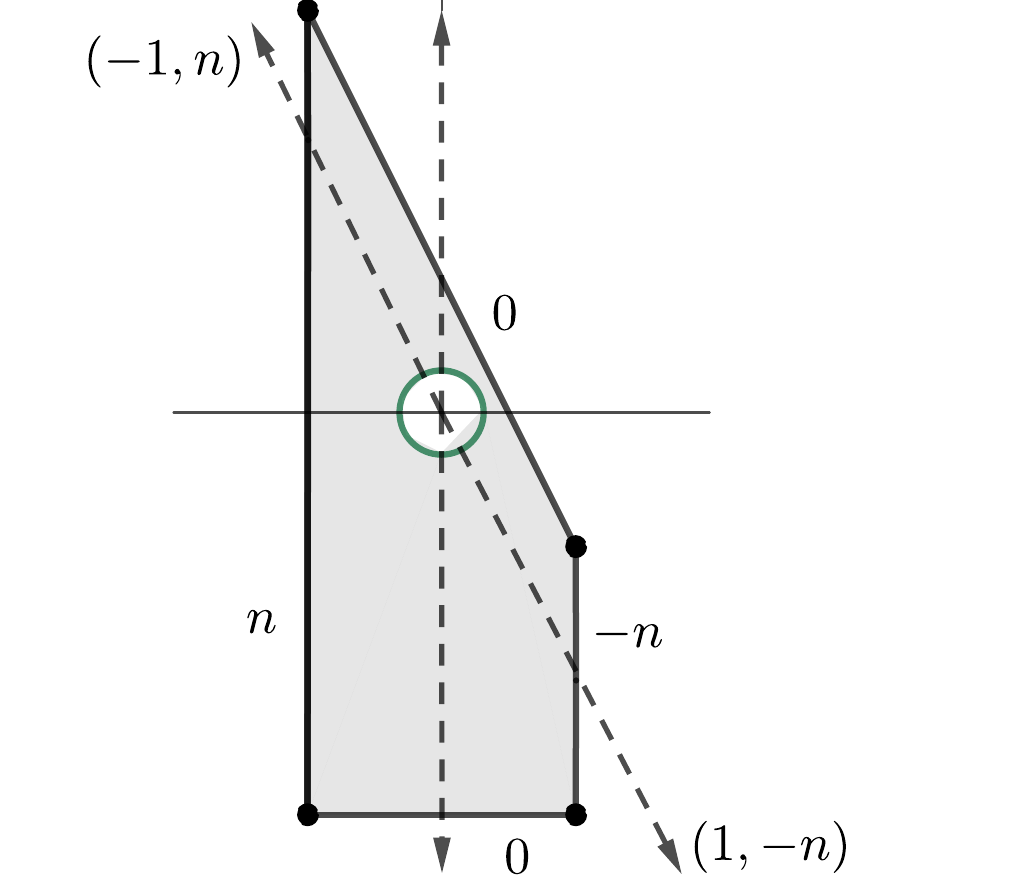}
   
    \end{minipage}
   \caption{A linear plumbing $(n,0,-n)$ and a cyclic plumbing $(n,0,-n,0)$}
\end{figure}

\begin{remark}\label{rem:forms} We make some observations concerning the symplectic structures on the given plumbings. The intersection form of the linear plumbing $(n,0,-n)$ is given by 
$$Q_n=\begin{bmatrix}
 n & 1 &0 \\
1 &  0 & 1 \\
0&1&-n
\end{bmatrix}.$$

 If $n\neq m$ and the parity of $n$ and $m$ is not the same or  $nm=0$ then the intersection forms $Q_n$ and $Q_m$ are not congruent, i.e. there does not exist a matrix $P\in GL(2, \mathbb Z)$ such that $Q_n=PQ_mP^T.$ Therefore, in that case, the corresponding linear plumbings are not even homeomorphic. On the other hand, for certain values of $n$ and $m$ the corresponding plumbing may be symplectomorphic. We do not proceed this elaboration, since, in this note, we focus on toric properties.
  \end{remark}

\noindent \textbf{2. Case} If $t_2-t_1<\pi $  the corresponding contact toric manifold is equivariantly contactomorphic to a lens space $L(k,l)$, $k>0, l\in \mathbb Z, $ with the unique universally tight contact structure $\xi_t$ and the moment cone shown on the left in Figure \ref{fig:lens}. As already shown in \cite[Lemma 5.1., Lemma 5.2.]{MS} this contact structure can be realised as a concave contact boundary of a linear plumbing $(s_1, \ldots, s_n)$, where

 \begin{equation} \label{eq:con_fraction}
 \frac{k}{l}=s_1-\frac{1}{s_2-\frac{1}{\cdots-\frac{1}{s_n}}},   \hskip2mm  \textrm{for some}  \hskip2mm s_1\geq0, s_2, \ldots,s_n\leq-2.
  \end{equation}

Consider now the contact toric manifold $(L(k,mk+l), \xi_t ),$  for any $m \in \mathbb Z.$ Similalry, it
 can be realised  as a concave contact boundary of 
 $(s^m_1, \ldots, s^m_{n_m})$, where

 \begin{equation} \label{eq:con_fraction2}
 \frac{k}{km+l}=s^m_1-\frac{1}{s^m_2-\frac{1}{\cdots-\frac{1}{s^m_{n_m}}}},  \hskip2mm  \textrm{for some}  \hskip2mm s^m_1\geq0, s^m_2, \ldots,s^m_{n_m}\leq-2. 
  \end{equation}

 On the other hand,  $(L(k,l),  \xi_t)$ is equivariantly contactomorphic to $(L(k,mk+l), \xi_t )$, for any $m \in \mathbb Z,$  through the  transformation 
$\begin{bmatrix}
 1 & m \\
0 &  1
\end{bmatrix}\in SL(2,  \mathbb Z)$
of the corresponding moment cones.
As $SL(2,  \mathbb Z)$ transformations preserve the self-intersection numbers of the spheres corresponding to edges, we conclude that
$(L(k,l), \xi_t), k>0, l\neq0$ can be realised  as a concave contact boundary of the linear plumbings  
$$(s^m_1, \ldots, s^m_{n_m}),\hskip1mm \textrm{for all} \hskip1mm m \in \mathbb Z.$$

As in the previous case, the corresponding symplectic toric manifolds are not equivariantly symplectomorphic since the sequences of self-intersection numbers of base surfaces differ. The later is true because these sequences give continued fraction expansions for distinct fractions. And, similarly as in Remark \ref{rem:forms}, 
for some choices of $k, l, m$ the underlying manifolds are not even homeomorphic, while for some other choices they can be even symplectomorphic.

\vskip2mm

\noindent \textbf{3. Case} If $\pi<t_2-t_1$ the corresponding contact toric structure $ \xi$  is obtained from $\xi_t$ by performing the sequence of $K \geq1$ half-Lutz twists.
According to \cite{MS}, if $(L(k,l),  \xi_t)$ is realised as a concave contact boundary of a linear plumbing $(s_1, \ldots, s_n)$ then
$(L(k,l),  \xi)$ can be realised as a concave contact boundary of a linear plumbing $(s_1, \ldots, s_n,\underset{2K}{\underbrace{0, \ldots, 0}})$. Since there are infinitely many plumbings  $(s_1, \ldots, s_n)$ for $(L(k,l),  \xi_t)$, we obtain  infinitely many plumbings for $(L(k,l),  \xi)$, keeping the number $K$ fixed.

\subsection{Contact  3-manifolds with a free toric action}
The proof relies on the previous case when $t_2-t_1=2N\pi$, for some $N \geq1$. There are infinitely many linear plumbings that are concave symplectic fillings of this contact toric structure, namely  
$$(n,0,-n,\underset{4N-2}{\underbrace{0, \ldots, 0}}), \hskip1mm \textrm{for all}\hskip1mm n \geq0.$$
To all these linear plumbings, we  perform the plumbing of the last vertex to the first vertex as explained in the proof of Proposition \ref{thm0}. That way we obtain infinitely  many cyclic plumbings whose boundary is equivariantly contactomorphic to $(T^3, \xi_N).$

\end{document}